\documentclass[12pt]{amsart}
\usepackage{amssymb}
\usepackage{amsfonts}
\usepackage{amsmath}

\usepackage{amsthm}

\newtheorem{theorem}{Theorem}[section]

\newtheorem{claim}[theorem]{Claim}

\newtheorem{corollary}[theorem]{Corollary}

\newtheorem{definition}[theorem]{Definition}

\newtheorem{lemma}[theorem]{Lemma}

\newtheorem{problem}[theorem]{Problem}
\newtheorem{proposition}[theorem]{Proposition}

\begin{document}

\title[Mackey countability]{On Docile Spaces, Mackey First Countable Spaces, and Sequentially Mackey First Countable Spaces}
\author[Bosch et al]{Carlos Bosch Giral     \and C\'{e}sar L.  Garc\'{i}a  \and   Thomas E. Gilsdorf \and Claudia G\'{o}mez - Wulschner \and Rigoberto Vera\\
Departamento de Matem\'{a}ticas, ITAM
R\'{\i}o Hondo \#1, Col. Progreso Tizap\'{a}n, M\'{e}xico, DF, 01080, M\'{e}xico}
\date{\today}

\maketitle

Email: \textbf{bosch@itam.mx, clgarcia@itam.mx, thomas.gilsdorf@itam.mx,
 claudiag@itam.mx, rveram@itam.mx}.

\begin{abstract}
\noindent In this article we discuss the relationship between three types of locally convex spaces: docile spaces, Mackey first countable spaces, and sequentially Mackey first countable spaces. More precisely, we show that docile spaces are sequentially Mackey first countable. We also show the existence of sequentially Mackey first countable spaces that are not Mackey first countable, and we characterize Mackey first countable spaces in terms of normability of certain inductive limits.   
\end{abstract}

 

\vspace{0.7in}

\section{Introduction.}     
\noindent   In this paper, {\it space} refers to   a Hausdorff locally convex topological vector space   $E = (E, \tau)$   over $\mathbb{K} = \mathbb{C} \, or \, \mathbb{R}$.   We assume that  $\,E$ is infinite dimensional unless otherwise indicated.    Unspecified notation follows standard locally convex spaces  notation such as in [2].  To denote $E$ as a linear subspace of another linear space $F$, we use $E \leq F$, and $E\lneq   F$ if $E$ is a proper subspace of $F$.       \\

\noindent The concept of a docile space was introduced in  [3]:   

\begin{definition}
A space    $E = (E, \tau)$    is {\bf  docile}  if every linear subspace of infinite dimension contains a  $\tau$- bounded set $B$ such that $span\{B\}$ is infinite - dimensional.        \\
\end{definition}

\noindent   In one of his seminal papers, G. W. Mackey [5]  defined the following:  

\begin{definition}

A space   $E = (E, \tau)$    is  {\bf first countable}, henceforth referred to as {\bf Mackey first countable},    if   for every sequence of $\tau-$ bounded sets  $(B_{n})$, there exists  a sequence of nonzero scalars $(a_{n})$,   such that  $\bigcup_{n=1}^{\infty } a_{n}B_{n}$ remains bounded.    \\
\end{definition}

 \noindent In  [1], the following was introduced  and the authors  proved that it is equivalent to the concept of a docile space:          \\

\begin{definition}

A space   $E$   is  {\bf sequentially Mackey first countable},  if given any sequence $\left( B_{n}\right) $ of bounded sets in $E$, and any
sequence $\left( x_{k}\right) $ from $\bigcup_{n=1}^{\infty }B_{n}$,  there exists a sequence $\left( \alpha _{k}\right) $ of scalars with infinitely many nonzero terms such that $\left\{ \alpha _{k}x_{k}:k\in  \mathbb{N}  \right\} $ is bounded.    \\

\end{definition}

\noindent  Observe that the definition is unchanged if one uses  sequences $\left( x_{k}\right)$ without reference to a  sequence $(B_{n})$ of bounded sets.    \\

\noindent Clearly every Mackey first countable space is sequentially Mackey first countable.  In this paper we will prove that there exist spaces that are sequentially Mackey first countable but not  Mackey first countable.  We will also obtain a characterization of Mackey first countability in terms of increasing sequences of bounded sets. \\

\section{Sequentially  Mackey  first countable spaces that are not Mackey first countable.}          

 \noindent  The existence of sequentially  Mackey  first countable spaces that are not Mackey first countable  conclusion will be in Corollary 2.2  which is a consequence of the following main result:      \\

\begin{theorem}
 The weak dual of a  nonnormable metrizable space  of  algebraic dimension greater or equal to $\aleph_{0}$  is not Mackey first countable.        

\end{theorem}  

\noindent In [4; (5), 393]   it is proved  that every nonnormable metrizable space is Mackey first countable and that the strong dual of an infinite dimensional  metrizable space, $(E^{\prime}, \beta(E^{\prime}, E))$, is not Mackey first countable.  In fact, it is easy to see that any regular  inductive (an inductive limit is {\bf regular}  if every bounded set is contained in and bounded in some step subspace) cannot be Mackey first countable.  \\

\noindent  In this paper we will make use of a bounding cardinal, $\mathfrak{b}$, defined as  $\mathfrak{b}=\min \{|F| \colon    F\subset \omega^{\omega}\,$\, such that for each $\, f\in \omega^{\omega}\,$\, there exists \, $\,g\in F\,$ for which $g$ is not less than $f \}$, where $\omega$ represents the ordinal of the natural numbers, $\mathbb{N}$, and $|F|$ denotes the cardinality of $F$.     \\

\noindent  It is known from ZFC, that  $\,\aleph_{1}\leq \mathfrak{b}  \leq \mathfrak{c} $, where $ \mathfrak{c}  = card(\mathbb{R})$.  \\

\vspace{2mm}
\noindent From  [3; Thm. 7. 4]  we have   that if the space $E$ is of algebraic dimension less than $\mathfrak{b}$, then its weak dual, $(E^{\prime}, \sigma(E^{\prime}, E)    )$ is docile.   Via  [1] we will obtain:     \\

\begin{corollary}
If   $E = (E, \tau)$   is metrizable  and of  dimension\, $\,\aleph_{o}\,$, then  $(E^{\prime}, \sigma(E^{\prime}, E))$  is sequentially Mackey first countable, but not Mackey first countable.    
\end{corollary}

\noindent We recall that  for an absolutely convex subset $C$ of a space $E$, the linear subspace generated by $C$ is $E_{C} = span\{C\} = \bigcup_{n=1}^{\infty} nC$.  Some simple observations of this construction  are:  If $a>0$, then $E_{a\cdot C} = E_{C}$, and if $A\subset C$, then $E_{A}  \leq  E_{C}$.  \\

\noindent The proof of Theorem 2.1  requires two preliminary lemmas:  \\

\begin{lemma}
Suppose  $\,C_{1}\subset C_{2}   \subset \cdots$ is a sequence of absolutely convex subsets of a space $E$, and that $\,a_{1},\, a_{2}, \cdots$ is a sequence of positive numbers.  If $\,C=\bigcup_{n=1}^{\infty}a_{n}C_{n}\,$,\, then   $\,E_{C}=\bigcup_{n=1}^{\infty}E_{C_{n}}\,$.    \\

\end{lemma}

\noindent {\it Proof}:  Because $\,a_{n}C_{n}\subset E_{C_{n}}\,$,  we have $\,C\subset \bigcup_{n=1}^{\infty}E_{C_{n}}\,$, and hence,   $\,E_{C}\subset \bigcup_{n=1}^{\infty}E_{C_{n}}\,$.     Next, let  $\,x\in \bigcup_{n=1}^{\infty}  E_{C_{n}}$. Then  $x\in E_{C_{n}}$, for some $n\in \mathbb{N}$. It follows that for some $r>0$,  $x\in rC_{n}\,$, which implies that there exists $y\in C_{n}$ with $ x=ry$.  From this,  \\

$$ x=ry=\frac{r}{a_{n}}(a_{n}y)\in \frac{r}{a_{n}}a_{n}C_{n}\subset \frac{r}{a_{n}}C\subset E_{C}.  \; \; \;  \Box $$

\vspace{0.25in}

If $C$ is bounded, then we can define a norm on $E_{C}$ via the Minkowski sublinear functional of $C$:  

$$\rho_{C}(x)= \inf\{r>0 :  x\in r\cdot C\}.   $$

\vspace{0.15in}

\begin{lemma}

Under the same hypotheses of the previous lemma, suppose additionally that each $\,C_{n}$  is bounded in $E$.  Consider the increasing sequence of  linear subspaces $E_{1} \subset E_{2} \subset \cdots$, where $E_{n}=(E_{C_{n}},\rho_{C_{n}}), \, n\in \mathbb{N}$.   Suppose that $F =\bigcup_{n=1}^{\infty} E_{n}$ and that  $C= \bigcup_{n=1}^{\infty}a_{n}C_{n}$  is bounded.  Then    $F =E_{C}$  and the topology induced by the norm   $\,\rho_{C}\,$    is weaker than the inductive limit topology formed by the sequence $( E_{C_{n}},\rho_{C_{n}})$.                               \\

\end{lemma}    

\noindent {\it Proof}: \/  The equality $F =E_{C}$ was proved in the previous lemma.  Now, observe that  $a_{n} C_{n} \subset a_{n} C_{n+1} \subset a_{n}C$, and that from  

$$E_{n}  = \bigcup_{i = 1}^{\infty}  i C_{n}.  $$



\noindent For each  $x\in E_{n}=\bigcup_{r>0} rC_{n}$,  $x=rz$  for some   $r>0$ and  $z\in C_{n}$.  From this, we have $\rho_{C_{n}}(x)\leq r$,  and  

$$    x=\frac{r}{a_{n}}(a_{n}z)\in \frac{r}{a_{n}}a_{n}C_{n}\subset \frac{r}{a_{n}}C. $$

\noindent  This last line implies that   $\rho_{C}(x)\leq \frac{r}{a_{n}}$, that is,   $a_{n}\rho_{C}(x)\leq r$.  By taking the infimum, we conclude that   $a_{n}\rho_{C}(x)\leq \rho_{C_{n}}(x)$ implies $ \rho_{C}(x)\leq \frac{1}{a_{n}}\rho_{C_{n}}(x)$, and we deduce  that  the topology induced by  $\rho_{C}$ on the subspace $E_{n}$  is weaker than the normed topology of  $\rho_{C_{n}}$, that is,   $\rho_{C_{n}}$ is weaker than   $\tau_{ind}$, which in turn gives us that $\rho_{C}$ is weaker than   $\tau_{ind}$.   \/  \/  $\Box$   \\

\noindent  {\it Proof of Theorem 2.1}.  \/  Let $\,V_{1}\supset V_{2}\supset \cdots$ be a decreasing countable base of zero neighborhoods of the space $E$.    We suppose that these neighborhoods are all absolutely convex and  $\tau$ - closed.  For each $n \in \mathbb{N}$, let $B_{n} = V_{n}^{\circ}$, the polar of the neighborhood $V_{n}$. By the Alaoglu - Bourbaki theorem, each $B_{n}$ is  $\sigma(E^{\prime}, E)$ - compact, hence, by  [4; (5.1),  262],  each $B_{n}$ is both bounded and complete with respect to the strong topology  $\beta(E^{\prime}, E)$.        \\

In  [4; (8), 394]   it is shown that   when $E$ is a nonnormable metrizable space,   $(B_{n})$ forms a fundamental sequence of bounded sets  in  $\beta(E^{\prime}, E)$, that for each $n\in \mathbb{N}$, $B_{n}$ does not absorb $B_{n+1}$, and that the linear subspace generated by each $\,B_{n}\,$, ($\,E'_{B_{n}}\,$) is a proper subspace of that generated by $B_{n+1}$.       \\ 

Of course,  $\bigcap_{n=1}^{\infty}V_{n}=\{0\}$  in $E$ implies that  $\,E^{\prime} =\bigcup_{n=1}^{\infty}B_{n}=\bigcup_{n=1}^{\infty}E'_{B_{n}}\,$.        \\

From the previous statements we conclude that the sequence $(B_{n})$  indicates that  $(E^{\prime}, \beta(E^{\prime}, E))$ is neither Mackey first countable nor sequentially Mackey first countable.   We will prove that:  

\begin{claim}   $(E^{\prime}, \sigma(E^{\prime}, E))$ is not Mackey first countable.    \end{claim}

\noindent  {\it Proof of Claim}.  \/  Suppose  there exists a sequence   $(a_{n})$ of positive numbers  such that $\bigcup_{n=1}^{\infty}a_{n} B_{n}\,$\, is\, $\sigma(E^{\prime}, E)$ - bounded.  Consider 

$$B = \overline{ absconv \left\{ \bigcup_{n=1}^{\infty} a_{n}B_{n})\right\}}^{\, \sigma(E^{\prime}, E)},  $$

\noindent  the absolutely convex $\sigma(E^{\prime}, E)$ - hull of $\bigcup_{n=1}^{\infty}a_{n}B_{n}$.   We will show that the assumption that  $B$ is   $\sigma(E^{\prime}, E)$ - bounded leads to the conclusion that $E$ is finite dimensional.    \\

\noindent  By Lemma 2.3, $\,E'=\bigcup_{n=1}^{\infty}E'_{B_{n}}=E'_{B}\,$, which means that $B$ is absorbing in $E^{\prime}$.  Lemma 2.4  implies that the topology induced by the norm  $\,\rho_{B}\,$\, in\, $\,E'\,$ is weaker than $\,\tau_{ind}\,$.   Moreover, because $B$ is a  $\,\rho_{B}$ - zero neighborhood, there exists a $\tau_{ind}$ -  zero neighborhood $U$ such that $\,U\subset B\,$.  Thus, $B$ is a $\tau_{ind}$ - bornivore, in particular, $B$ is absorbing in $E^{\prime}$  \\

\noindent    On the other hand, the  $\sigma(E^{\prime}, E)$ - completeness of each $B_{n}$  implies that  each  $\,(E'_{B_{n}},\rho_{B_{n}})\,$ is a Banach space, giving us that the inductive limit of these spaces is barrelled (cf [2;  214]).  Thus, as $\,B\subset E'\,$ is an absolutely convex, absorbing, $\sigma(E^{\prime}, E)$ - closed  (hence, $\,\tau_{ind}$- closed) set, i.e., a barrel, $B$ is a zero neighborhood with respect to  $\,\tau_{ind}$.  We conclude that $\beta(E^{\prime}, E)   =\tau_{ind}\,$\, in  $\,E'\,$.    Meanwhile,    ($\,E''=(E'[\beta'])'$  is metrizable: [4; (5),  301] and each $\sigma(E'', E)$ - bounded subset  is $\,\beta'$-equicontinuous ([2; 3.6.2,  212]).  Because $U$  is absorbing in $E'$, its polar in  $\,(E', \tau_{ind}')^{\prime} =E''\,$ is $\sigma(E'', E)$ - compact, which means $U^{\circ}$ is $\beta(E'', E)$ -  bounded.  From $U\subset B$,   $B^{\circ}   \subset U^{\circ}\subset E''$, implying that $B^{\circ}$ is  $\beta(E'', E)$ - bounded, and  $\sigma(E'', E)$ - bounded as well.  Because  $B$ is    $\sigma(E^{\prime}, E)$ - bounded,   $B^{ \circ} \subset E^{\prime \prime}$ is absorbing, that is, $\,E''_{B^{\circ}}=E''\,$.  Furthermore,  $\,B^{\circ} \subset E''\,$ is $\,\beta'$-equicontinuous and thus  $\,C=E\cap B^{\circ}\subset E\,$\, is\, $\,\beta'$-equicontinuous, hence, $\tau$ -  precompact.  The equicontinuity of $C$ gives us that $C$ is a $\tau$ -  zero neighborhood, and therefore, $E$ is finite dimensional.      This contradiction completes the proof.   \/ \/  $\Box$   \\


\noindent {\it Proof of Corollary 2.2}.   This follows from the Theorem 2.1  above, and from the proof that such a space is docile ([3; Thm. 7])  given that  $\,\aleph_{o}< \mathfrak{b}$.   \/ \/   $\Box$   \\

\noindent   In  [1] the authors  proved  that a Hausdorff quotient of a Mackey first countable space is sequentially Mackey first countable.   Theorem 2.1  leaves that conclusion as a partial answer to G. W. Mackey's original question, which, in its original form, remains open:   \\    

\begin{problem}  (G.W. Mackey,  [5]).  Is a Hausdorff quotient of a Mackey first countable locally convex space again Mackey first countable?
\end{problem}

\section{A characterization of Mackey first countable spaces.}

\noindent  Let $\,B_{1}\subset B_{2}\subset , \cdots$ be an increasing sequence of absolutely convex closed bounded sets in a space $E$.  Put  $\,(F=\bigcup_{n=1}^{\infty}E_{B_{n}},\,\tau_{ind})\,$\, as the inductive limit of the normed  linear subspaces generated by this sequence. By [2; 7.3.4, 222], the space $\,(F,\tau_{ind})\,$ is bornological.  In the main  result of this section  we will prove that $E$ is Mackey first countable if and only if $\,(F,\tau_{ind})\,$ is normed.     \\




\begin{proposition}
Suppose $\,\{B_{n}\}_{n\in \Bbb{N}}$ is any sequence of closed, convex subsets of $E$, such that $B_{1}$ is balanced.  Define, for each  $\,n\in \Bbb{N}$,  $\,B_{n}^{\#}  =\bigcup_{k=1}^{n}B_{k}$.    If  $\,\{a_{n}\}_{n\in \Bbb{N}}$ is such that  $0 < a_{n} < 1$ for each $n \in \mathbb{N}$, then  $\,\bigcup_{n=1}^{\infty}a_{n}B_{n}\,$ is bounded if and only if, $\,\bigcup_{n=1}^{\infty}a_{n}B_{n}^{\#} \,$ is bounded.   

\end{proposition}

\noindent {\it Proof}.   ($\,\Leftarrow\,$):    Because  $\,B_{n}\subset B_{n}^{\#}$, this implication is clear.           \\

\noindent   ($\,\Rightarrow\,$):  Without loss of generality, we may suppose that  $a_{1}> a_{2}> \cdots $.   Let  $\,A\subset E\,$ be $\,\tau$ - bounded such that  $\,\bigcup_{n=1}^{\infty}a_{n}B_{n}\subset A\,$. Observe that    \\

$$\,a_{1}B_{1}=a_{1}B_{1}^{\#}, \;  a_{2}B_{2}^{\#} =a_{2}(B_{1}\cup B_{2})= a_{2}(B_{1}\cup a_{2}B_{2})\subset a_{1}B_{1}\cup a_{2}B_{2}\subset A.  $$

\noindent Thus, 
$$\,a_{n}B_{n}^{\#}  =a_{n}(\bigcup_{k=1}^{n}B_{k})=\bigcup_{k=1}^{n}a_{n}B_{k}\subset \bigcup_{k=1}^{n}a_{k}B_{k}\subset A,    $$

\noindent which shows that  $\,\bigcup_{k=1}^{\infty}a_{k}B_{k}^{\#}    \subset A$. \/ \/  $\Box$ \\

\begin{corollary}
Suppose   $\,\{B_{n}\}_{n\in \Bbb{N}}$ is any sequence of  closed, convex subsets of $E$, with  $B_{1}$  balanced.  Define  $B_{n}^{\#}$ by:

$$  \,B_{n}^{\#}  =\overline{absconv \left\{ \bigcup_{k=1}^{n}B_{k} \right\}}^{\, \tau}.$$

\noindent   If  $\,\{a_{n}\}_{n\in\Bbb{N}}$ is such that  $0 < a_{n} < 1$ for each $n \in \mathbb{N}$, then  $\,\bigcup_{n=1}^{\infty}a_{n}B_{n}\,$ is bounded if and only if, $\,\bigcup_{n=1}^{\infty}a_{n}B^{\#}_{n}\,$ is bounded.   

\end{corollary}

\noindent {\it Proof}.   If $\,\bigcup_{n=1}^{\infty}a_{n}B^{\#}_{n}\,$ is bounded, then  we obtain again,\, $\,B_{n}\subset B_{n}^{\#}$.     \\

\noindent Suppose  $\,\bigcup_{n=1}^{\infty}a_{n}B_{n}\,$ is bounded.   Consider  $\,A\subset E$, an absolutely convex closed, bounded set, such that $\,\bigcup_{n=1}^{\infty}a_{n}B_{n}\subset A\,$.  By the previous proposition,   \\

$$  \bigcup_{n=1}^{\infty}a_{n}(B_{1}\cup\cdots  \cup B_{n})\subset A.   $$

\noindent  We observe that  $$ \overline{absconv \left\{  a_{n}(B_{1}\cup\cdots \cup B_{n}) \right\}}^{\, \tau}   =a_{n}\,B_{n}^{\#}.   $$

\noindent  It follows  that  $\,a_{n} B_{n}^{\#}   \subset A\,$, for each $\,n\in \Bbb{N}\,$, hence, 

$$\bigcup_{n=1}^{\infty}a_{n}B_{n}^{\#}   \subset A.  \; \; \Box  $$

\vspace{0.25in}
\begin{corollary}
$$ \sum_{n\in \Bbb{N}}E_{B_{n}}=\bigcup_{n=1}^{\infty}E_{B_{n}^{\#}   } \leq  E.  $$
\end{corollary}

\noindent {\it Proof}.   For each  $\,n\in \Bbb{N}\,$ we have  $\,E_{B_{n}} \leq  E_{B_{n}^{\#}  }$, which gives us 

$$\,\sum_{n\in \Bbb{N}}E_{B_{n}} \leq  \bigcup_{n=1}^{\infty}E_{B_{n}^{\#}}.    $$

\noindent  Let   $\,x\in \bigcup_{n=1}^{\infty}E_{B_{n}^{\#}}$. Then   $\,x\in E_{B_{n}^{\#}}\,$ for some  $\,n\in \Bbb{N}\,$.  By 
[4; (1), 173], $x=\sum_{k=1}^{m}\alpha_{k}x_{k}\,$  with $\,\alpha_{k}\geq 0\,$,   $\,\sum_{k=1}^{m}\alpha_{k}=1\,$\, and  $\,x_{k}\in \bigcup_{k=1}^{m}E_{B_{k}}\,$ for each 
$\,k=1,2,\cdots  ,m\,$.  We obtain $\,x\in \sum_{n\in \Bbb{N}}E_{B_{n}}\,$.  \/ \/ $\Box$  \\

\begin{theorem}
The space  $(E, \tau)$ is Mackey first countable if and only if, for each increasing collection  $\,\{B_{n}\}_{n\in \Bbb{N}}$   of bounded, absolutely convex closed sets in $E$, the corresponding inductive limit,  $\,ind_{n\in \Bbb{N}}(E_{B_{n}^{\#}},\rho_{B_{n}^{\#}})$, is normable.   
\end{theorem}

\noindent {\it Proof}.   Let $\,F=\bigcup_{n=1}^{\infty}E_{B_{n}^{\#}} \leq  E\,$, and let $\,\tau_{ind}\,$ denote the corresponding inductive limit topology on $F$.    \\

\noindent  Assume that for each increasing collection  $\,\{B_{n}\}_{n\in \Bbb{N}}$   of bounded, absolutely convex closed sets in $E$,  $\,ind_{n\in \Bbb{N}}(E_{B_{n}^{\#}},\rho_{B_{n}^{\#}})$, is normable.  Let  $\,\{B_{n}\}_{n\in \Bbb{N}}$ be a collection of bounded subsets of $E$.  Using $\,\rho\,$ to denote the normed topology on $F$, we have $\,\tau_{ind}=\rho\,$.  If $\,D\subset F\,$ is the $\rho$ - closed unit ball, the $(\forall n \in \mathbb{N})(\exists a_{n} >0)$ such that $\,a_{n}B_{n}^{\#}   \subset D\cap E_{B_{n}^{\#}}\subset D\,$.   From this,  $\,\bigcup_{n=1}^{\infty}a_{n}B_{n}^{\#}   \subset D\,$.   Moreover, $\,B_{n}\subset B_{n}^{\#}$ implies that  $\,\bigcup_{n=1}^{\infty}a_{n}B_{n}\subset D\,$.   Meanwhile, $D$ is $\rho$ -  bounded,  with $\,\tau \leq \tau_{ind}=\rho\,$, and this tells us that $D$ is $\tau$ -  bounded as well.   \\

\noindent  Conversely, if   $(E, \tau)$ is Mackey first countable, then we will prove that there exists a norm $\rho$ on $F$ such that $\,\tau_{ind}=\rho\,$.   To this end, let $\,B_{1}\subset B_{2}\subset \cdots$ be an increasing sequence of closed, absolutely convex bounded subsets of $E$.  Then there exists a sequence  $\,\{a_{n}\}_{\Bbb{N}}$ in  $(0,1]$, with $\,a_{1}> a_{2}> \cdots > 0\,$ such that  $\,\bigcup_{n=1}^{\infty}a_{n}B_{n}$ is bounded in $E$.  Put   

$$ B =  \overline{absconv\left\{  \bigcup_{n=1}^{\infty}a_{n}B_{n} \right\}}^{\, \tau}.  $$

\begin{claim}
$\,\rho_{B}=\tau_{ind}\,$.
\end{claim}

\noindent {\it Proof of claim}.   We have that $B$ is absolutely convex, closed and bounded.  By Lemma 2.4,   $\,F=E_{B}\,$.  On the other hand, by the definition of the inductive limit topology, $\,\rho_{B}\leq \tau_{ind}\,$.  Hence,   $\,(F,\rho_{B})\,$ and $\,(F,\tau_{ind})\,$ are both bornological.  We conclude via [2; 3.7.1 (a), 220], that     $\,\rho_{B}=\tau_{ind}\,$.   \/ \/ $\Box$  \\

\begin{proposition} The previous result applies to any  countable collection of absolutely convex closed and bounded subsets.  \end{proposition}

\noindent  {\it Proof}.    Consider the family $\wp$  of finite subsets of $\mathbb{N}$.  The family $\wp$ is ordered by inclusion.  Moreover,  given $\,S\in \wp\,$, $\,C_{S}=\sum_{k\in S}E_{B_{k}}\,$ satisfies:  If $\,S\subset S'\,$, then    $\,C_{S} \leq  C_{S'}\,$.  Therefore, we can speak of the inductive limits  of  such linear subspaces of the form   $\,\{C_{S}\}_{S\in \wp}\,$:  $\,G=ind\, C_{S}\,$, with the respective topology $\,t_{ind}\,$.     \\    

\noindent  By Corollary 3.3, $F = G$, and because $\,C_{S}\,$ is a linear subspace of some  $\,E_{B_{n}^{\#}}\,$, we have:  $\,\tau_{ind}\leq t_{ind}\,$.  \/ \/ $\Box$          \\

  \vspace{0.5in}




\underline{ACKNOWLEDGEMENT}: \ The authors  would like to acknowledge
research support for this paper  by the Asociaci\'{o}n Mexicana de Cultura, A. C.      \\

\bigskip 

\begin{center}
\underline{\textbf{REFERENCES}}
\end{center}

\bigskip

\noindent  [1]  Bosch Giral, C.,  Gilsdorf, T.,  G{\'o}mez-Wulschner, C.  {\it Mackey first countability and docile locally convex spaces}.   Acta Math. Sin. (Engl. Ser.) {\bf 27},  (2011), no.4,  737 - 740. \\

\noindent   [2]  Horv\'{a}th, J. \underline{Topological Vector Spaces and Distributions, Vol. I}.  Addison - Wesley, (1966).\\

\noindent  [3]   Kakol, J.,  Saxon, S. A.,   Todd, A. R.  {\it Docile locally convex spaces}, Contemp. Math., \textbf{341}, AMS, (2004), 73 - 77.       \\

\noindent [4]   K\"{o}the, G.  \underline{Topological Vector Spaces I}.  Springer Verlag, NY, (1969).      \\

 \noindent   [5]  Mackey, G. W.  \textit{On infinite- dimensional linear spaces}. Trans. AMS, \textbf{66}, no. 57 (1945), 155 - 207.    \\

\noindent [6]  P\'{e}rez Carreras, P., Bonet, J. \underline{Barrelled Locally Convex Spaces}. \ North Holland Math. Studies, {\bf 131}, (1987).    \\
  
\noindent  [7]  Saxon, S. A., S\'{a}nchez Ruiz, L. M.  {\it Optimal cardinals for metrizable barrelled spaces}.  J. London Math. Soc. (2),  {\bf 51},  (1995),  no. 1, 137 - 147.           \\

\end{document}